\documentstyle [12pt]{article}
\newcommand{\sect}[1]{\setcounter{equation}{0}\section{#1}}

\newcommand{\f}{\frac}

\newcommand{\p}{\partial}
\newcommand{\FR}{Fun(\rn_q^N)}

\newcommand{\e}{\vec{e}}
\newcommand{\y}{\vec{y}}
\newcommand{\w}{\vec{w}}

\newcommand{\J}{\vec{j}}
\newcommand{\Pg}{\vec{\pi}}
\newcommand{\ot}{\otimes}
\newcommand{\La}{\Lambda}

\newcommand{\ve}{{\varepsilon}}
\newcommand{\ap}{\approx}
\newcommand{\bc}{\begin{center}}
\newcommand{\ec}{\end{center}}
\newcommand{\be}{\begin{equation}}
\newcommand{\ee}{\end{equation}}

\newcommand{\M}{{\bf L}}
\newcommand{\k}{{\bf k}}

\def\lcross{{>\!\!\!\triangleleft}}

\newcommand{\cn}{{\bf C}}
\newcommand{\rn}{{\rm\bf R}}
\newcommand{\zn}{{\bf Z}}
\newcommand{\nn}{{\rm\bf N}}
\newtheorem{prop}{Proposition}

\newtheorem{theorem}{Theorem}

\begin{document}

{\huge{\bf The q-Euclidean algebra $U_q(e^N)$
              and the corresponding q-Euclidean lattice}}

{}~
\bc
{\bf Gaetano Fiore\footnote{Alexander-von-Humboldt fellow}}
\ec
\bc
{\it Sektion Physik der Universit\"at M\"unchen,}
\ec
\bc
{\it Theresienstrasse 37, D-80333 M\"unchen, Germany}
\ec
\bc
LMU-TPW 95-4, February 1995,
\ec
\section*{\center Abstract}

We present the Euclidean Hopf algebra $U_q(e^N)$ dual of $Fun(\rn_q^N\lcross
SO_{q^{-1}}(N))$
and describe its fundamental Hilbert space representations \cite{fioeu}, which
turn out
to be rather simple `` lattice-regularized '' versions of the classical ones,
in the sense that the spectra of squared momentum components are discrete and
the corresponding eigenfunctions normalizable.
A suitable notion of classical limit is introduced, so that we recover the
classical
continuous spectra and generalized (non-normalizable) eigenfunctions in that
limit.

\vskip1truecm

\section*{Introduction}

{}~~~
Since their birth quantum groups \cite{dr} have found a number of different
applications to physics and mathematics. In particular
they can be used to generalize the ordinary notion
of space(time) symmetry. This generalization is tightly coupled to a
radical modification of the ordinary notion of space(time) itself.
{}From this viewpoint inhomogenous group symmetries such as Poincar\'e's and
the
Euclidean one yield physically relevant candidates for quantum group
generalizations; Minkowski space $M^4$ and Euclidean
$\rn^N$ one are then the corresponding space(time) manifolds.
One can  generalize the latter by
the $N$-dimensional ($N\ge 3$) Euclidean space $\rn_q^N$
\cite {frt}, its symmetry by the q-Euclidean one carried by the
Hopf-algebra $E_q^N:=\rn_q^N\lcross SO_q(N)$ $\cite{maj2,schl,maj3}$ or
equivalently by its dual \cite{maj2,fioeu}, which here we will call $U_q(e^N)$.
In Ref \cite{fioeu} we classified the fundamental Hilbert space representations
of $U_q(e^N)$;
here we represent the latter results in a more pedagogical and explicit way and
add some new ones.

A major physical motivations for such generalizations is
the desire to discretize space(time) (or momentum space) in a `` wise ''
way for QFT regularization purposes. Nowadays such a discretization is usually
performed by
approximating the points of the space(time) (or momentum space) continuum by
the
points of a
lattice. In the case of the cubic Euclidean lattice, for instance, the
coordinates $x^i$
($i=1,2...,N$) can assume only the values $an^i$, where $a$ is the lattice
spacing and
$n^i\in \zn$;
one chooses as a basis of the Hilbert space ${\cal H}$ of physical states the
set
$\{|n^1,...,n^N>\}_{n_i\in\zn}$ of eigenvectors of the $N$ commuting
observables
$x^i$
with eigenvalues $an^i$. On the other hand, it is known that standard
lattices used in regularizing QFT do not carry representations of
discretized versions (in the form of discrete subgroups) of the
associated inhomogenous groups; actually, the notion of a group is too tight
for this scope. For instance, the Euclidean
cubic lattice is invariant only under a discretized version of the translation
subgroup of the Euclidean group,
but not of the rotation one; in other words, we are able only to represent the
latter subgroup
on ${\cal H}$. On the contrary, the notion
of symmetry provided by quantum groups is broad enough
to allow the existence of lattices whose points are mapped into each other
under the action of the whole inhomogeneous q-groups. The main purpose of this
paper
is to describe how this occurs in the case of the q-Euclidean symmetry and how
in the
limit $q\rightarrow 1$ one recovers the ordinary representation spaces. One
concludes
that the $q$-Euclidean lattice introduced in ref. \cite{fioeu} seems very
appealing in view of full covariant
regularizations of Euclidean QFT; actually,
a $q$-deformed version of the $\varepsilon$-tensor on $\rn_q^N$ is also
available
\cite{fio3,fioeu,maj4},
so allowing the construction of the pseudo-tensors which are needed for chiral
field theories.

The main difference w.r.t  the cubic lattice stems from the following fact.
The $N$ configuration-space coordinates $x^i$
(as well as the momenta $p^i$) don't commute
with each other; therefore we can use a complete set of commuting observables
consisting only
partially (i.e. for about one half) of functions of the $p^i$ (or,
alternatively, of the $x^i$) and,
as for the rest, by angular momentum components. Their spectra are discrete.
The
lattice in
the present situation has namely $[\f {N+1}{2}]$ dimensions in $p$-space and
$[\f {N}{2}]$
dimensions ($[a]$ denotes the integer part of $a$)
in angular momentum space; to each point of the lattice there corresponds a
unique eigenvector
belonging to a basis of ${\cal H}$ and labeled by $N$ integers. Notably, under
the action of the
generators of the q-Euclidean algebra each vector is mapped simply into a new
one with labels
differing at most by $\pm 1$. In the sequel we will consider as algebra of
observables the one
generated by $p^i$'s and the angular momentum components, since they generate
the physically
relevant q-deformed Euclidean algebra (q-translations + q-rotations);
but both the commutation relations and the representation
theory would be exactly the same (under the replacement $x^i\rightarrow p^i$)
if we considered the $x^i$ instead.

In section 2 we briefly introduce the q-deformation $U_q(e^N)$
of the universal enveloping algebra of the Euclidean Lie algebra $e^N$ which we
are going
to adopt as quantum symmetry.
We will be quite explicit in the case $N=3,4$, for which we also write down the
analog of
the Pauli-Lubanski casimirs.

$U_q(e^N)$ is the Euclidean analogue
of the q-deformed Poincare' Hopf algebra (of u.e.a. type) Ref.
\cite{wess,maj2}.
In both cases the inhomogeneous Hopf algebra contains the
homogeneous one as a Hopf subalgebra which can be obtained from it by setting
$p^i=0,\La=1$ ($\La$ is the ``dilaton"), and all commutation relations
are homogeneous in $p$, contrary to what happens for inhomogenous
Hopf algebras obtained through contractions \cite{cel,luk1,luk2}.
Representation
theory is
also developed in a similar way as in ref. \cite{wess}.

Section 3 is devoted to a detailed description of fundamental (i.e. irreducible
one-particle) Hilbert space $*$-representations of $U_q(e^N)$
(we will call them `` irreps '' in the sequel). The case $N=3$ is analysed
first, as an
introduction to the general case. We choose a Cartan subalgebra
(i.e. a complete set of commuting observables) consisting basically of
two parts, $[\f {N+1}2]$ squared momentum components and $[\f N2]$
angular momentum components ($[a]$ denotes the
integer part of $a$).  The points of the spectra make up a q-lattice.
One important fact is that the irreps turn out to be
of highest weight type. Moreover, they can
be obtained from tensor products of the singlet one (i.e.
the one describing a particle with zero $U_q(so(N))$-highest weight) and some
representation
of $U_q(so(N))$; for instance, the irreps with $N=3$ are obtained from the
tensor
product of the q-boson (i.e. zero spin) representation of $U_q(e^3)$ with a
representation
of some spin $j\in\nn$ of $U_q(so(3))\ap U_q(su(2))$, in analogy with the
undeformed case.
The spectra of all observables are discrete, in particular the spectra of
squared momentum
components, as expected. The corresponding eigenvectors
are normalizable and make up an orthogonal basis of the Hilbert space
of each irrep.
A cumbersome `` kinematical PT (parity + time-inversion) asymmetry '' appears
in the structure of the spectra of the angular momentum observables;
it disappears in the limit $q\rightarrow 1^-$.

In section 4 we clarify in which sense the Euclidean algebra/representations go
to the
classical ones in the limit $q\rightarrow 1$. In the classical representation
we
know that
the eigenvectors of operators which are only functions of the momenta are
distributions,
tipically they are
 delta-functions in momentum space. We show how to construct $q$-dependent
integer labels $n_i(q)$
and coefficients $\alpha(q)$
such that $\alpha(q)|n_i(q)>_q$ (eigenvectors belonging to the
$q$-representation)
are delta-convergent functions in the limit $q\rightarrow 1$.

We can think of the irreps studied in section 3 as
describing the (time-independent) dynamics of
a free nonrelativistic particle with arbitrary `` generalized ''
$U_q(so(N))$-spin on $\rn_q^N$.
The subalgebra $\hat U_q(e^N):=U_q(e^N)/(\La-1)$
can be considered as the quantum group symmetry of the hamiltonian
\be
H:=\f{(p\cdot p)}{2M}.
\ee
of the system; therefore
all states with a given energy should be obtained from each other
by the action of $\hat U_q(e^N)$, as in the classical case,
different eigenspaces of the energy should be obtained from each other by the
action of the dilatation operators $\La^{\pm 1}$.

Some notational remarks are necessary before the beginning.
For representation purposes we will assume in section 3 that
$q\in {\bf R}^+$, and we will limit ourselves to the case $0<q\le 1$;
the case $q>1$ can be treated in an analogous way.
We set $h=h(N)=\cases{0~~~if~~N=2n+1\cr 1~~~ if~~N=2n\cr}$ to allow
a compact way of writing relations valid both for even and odd $N$.
Unless stated differently, in our notation a space index $i$ can
take all the integer values between $-n$ and $n$ including/excluding $i=0$
if $N=2n+1,2n$ respectively.
When $N=2n$ there is a complete invariance of the validity of all
the results
under the exchange of indices $i=-1\leftrightarrow i=1$, so that we will
normally omit writing down explicitly the results that can be obtaind
by such an exchange. We will often use the shorthand notation
$[A,B]_a:=AB-aBA$ ($\Rightarrow [\cdot,\cdot]_1=[\cdot,\cdot]$).
Indices are raised
and lowered through the q-deformed metric matrix $C:=||C_{ij}||$, for instance
\be
a_i=C_{ij}a^j,~~~~~~a^i=C^{ij}a_j,~~~~~~~~~~~~~~~~
C_{ij}:=q^{-\rho_i}\delta_{i,-j},
\ee
where
\be
(\rho_i):=\cases {
(n-\f{1}{ 2},n-\f{3}{ 2},...,\f{1}{ 2},0,-\f{1}{ 2}...,\f{1}{ 2}-n)
{}~~~~~~~~~~if~N=2n+1 \cr
(n-1,n-2,...,0,0,...,1-n)~~~~~~~~~~~~~~~if~N=2n. \cr}
\ee
$C$ is not symmetric and coincides with its inverse: $C^{-1}=C$.

\sect{The Euclidean $*$-algebra $U_q(e^N)$}

The Hopf algebra which we are going to use, $U_q(e^N)$, was constructed in Ref.
\cite{maj2}
and in equivalent form in ref. \cite{fioeu}  by an inhomogeneous extension of
the
Hopf algebra $U_q(so(N))$ of  ``infinitesimal q-rotations'' (in analogy with
the
undeformed
construction). $U_q(e^N)$ is  the Hopf dual of $Fun(\rn_q^N\lcross
SO_{q^{-1}}(N))$ \cite{schl, maj2}.
In  ref. \cite{fioeu} work, we added to the Drinfeld-Jimbo generators of the
latter
first the q-derivatives on $\rn_q^N$ as infinitesimal generators $p^i$ of
q-translations
and then one more generator $\La$, generating dilatations; the coalgebra and
antipode for $U_q(e^N)$
were derived from the Leibnitz-rule of q-differential operators
(the role of the coalgebra in representation theory is to allow the
construction
of many-particle
representations starting from one-particle ones). On the algebra $U_q(e^N)$
there exists a notion of complex conjugation $*$, (which will play the role of
hermitean
conjugation of operators).
However, since the coalgebra is uncompatible,
at least in the usual sense, with the $*$-structure, here we focus the
attention
on
the algebra structure of $U_q(e^N)$ which we need to develop the theory of
one-particle
representations.

\subsection{A Chevalley basis of $U_q(so(N))$}

{}~~~A Cartan-Weyl basis of $U_q(so(N))$
is the set $\{\M^{ij},(\k^l)^{\pm \f 12}\}$ ($i< j,\neq -j;~n\ge l\ge 1$) with
commutation relations given
below. Its elements were realized in Ref. \cite{fio4} as q-differential
operators on $\rn_q^N$;
this is the q-deformed analogue of realizing the generators of $so(N))$ as
``angular momentum components''.
To help the reader in the identification of the corresponding classical angular
momentum components,
we give here their classical limits
\be
\M^{ij}\stackrel{q\rightarrow 1}{\longrightarrow}{x^i\p^j-x^j\p^i},~~~~~~~~~~~~
\f{\k^l-1}{q^2-1}\stackrel{q\rightarrow
1}{\longrightarrow}{x^l\p^{-l}-x^{-l}\p^l},
\ee
where $x^i,\p^j$ denote the classical coordinates/derivatives,
$\p^ix^j=\delta^{i,-j}+x^j\p^i$; the latter are chosen not to be
real, but complex combinations such that $(x^i)^*=x^{-i}$, $(\p^i)^*=-\p^{-i}$.
According to this construction, $U_q(so(N))$ is realized as a subalgebra of the
differential
algebra on $\rn_q^N$.

The $\k^i$'s generate a Cartan subalgebra of $U_q(so(N))$.
The elements $\M^{-i,i+1},\M^{-i-1,i},\k^i(\k^{i+1})^{-1}$
(together with $\M^{12},\M^{-2-1},\k^1\k^2$ in the case $N=2n$) for
$i=h,h+1,...,n$ are `` Chevalley generators '' (i.e. algebraically independent
generators)
of $U_q(so(N))$ coinciding \cite{fio4} with the Drinfeld-Jimbo ones,
up to some rescaling of the roots $\M$ by suitable
functions of $\k^i$ ($\k^0\equiv 0$). The correspondence between the
Chevalley generators $\M^{-i,i+1}$ corresponding to positive roots and
the spots of the Dynkin diagram of $so(N)$ is shown in fig. 1.
All the other generators $\M^{ij}$ can be constructed starting from
them as follows:
\be
[\M^{-jl},\M^{-lk}]_q=q^{\rho_l}\M^{-j,k}~~~~~~~~~~[\M^{-kl},\M^{-l,j}]
_q=q^{\rho_l+1}\M^{-k,j},~~~~~~~~~~~~~~n\ge k>l>j\ge -h(N)
\ee
\be
[\M^{l-1,k},\M^{1-l,l}]_{q^{-1}}=q^{\rho_l-1}\M^{lk}~~~~~~~~~~
[\M^{-l,l-1},\M^{-k,1-l}]_{q^{-1}}=q^{\rho_l}\M^{-k,-l}~~~~~~~~~~~~~~~~~~~~~
2\le l<k\le n
\ee
\be
[\M^{0k},\M^{01}]=q^{-1}\M^{1k}~~~~~~~~~~[\M^{-10},\M^{-k0}]=\M^{-k,-1}
{}~~~~~~~~~~~~~~~~~~~~~~~1<k\le n~~~if~~N=2n+1;
\ee
these relations can be easily verified by the reader in the limit $q=1$ using
the limits (1.1).

Once introduced the basis $\{\M^{ij},\k^l\}$ ($i< j,\neq -j;~n\ge l\ge 1$),
then the commutation relations satisfied by the Chevalley generators can be
summarized in the following way.
\begin{itemize}
\item Commutation relations between the generators of the Cartan
 subalgebra and
the simple roots:
\be
[\k^i,\M^{\pm (1-k),\pm k}]_a=0~~~~~~a=\cases{q^{\pm 2}~~if~i=k\le n \cr
q^{\mp 2}~~~if~~i=k-1 \cr 1~~~otherwise}~~~~~~~~~~~~~~~~~~
[\k^i,\k^j]=0;
\ee
\item commutation relations between positive simple roots
(the ones appearing on the left of the q-commutators) and negative ones ( the
ones
appearing on the right):
\be
[\M^{1-m,m},\M^{-k,k-1}]_a=0~~~~~~~~~~~~~~~a=\cases{q^{-1}~~~~m\pm 1=k \cr
1~~~~if~~k\neq m,m\pm 1\cr}~~~~m,k\ge h(N)+1,
\ee
\be
[\M^{12},\M^{-2,1}]=0~~~~~~~~~~[\M^{-1,2},\M^{-2,-1}]=0~~~~~~~~~~~
{}~~~~~~~~~~if~~~~~N=2n,
\ee
\be
\cases{[\M^{1-m,m},\M^{-m,m-1}]_{q^2}=q^{1+2\rho_m}\f{1-\k^{m-1}(\k^m)^{-1}}
{q-q^{-1}}~~~~~~~~2\le m\le n\cr
[\M^{01},\M^{-1,0}]_q=q^{-\f 12}\f{1-(\k^1)^{-1}}{q-q^{-1}}~~~~~~~~
if~~N=2n+1;\cr}
\ee
\item Serre relations:
\be
[\M^{1-m,m},\M^{1-k,k}]=0~~~~~~~~~[\M^{-m,m-1},\M^{-k,k-1}]=0~~~~~~~~~~~~~~~
m,k>0,~~~~ |m-k|>1
\ee
\be
[\M^{1-k,k} ,\M^{2-m,m}]_a=0=[\M^{-m,m-2},\M^{-k,k-1}]_a
{}~~~~~~~~~~~~~~a=\cases{q~~~if~~k=m\cr q^{-1}~~~
if~~k=m-1\cr}~~~~~~~~m\ge 3
\ee
\be
\cases{[\M^{01},\M^{12}]_{q^{-1}}=0 \cr [\M^{-1,2},\M^{02}]_q=0\cr}
{}~~~~~~~~~
\cases{[\M^{-2,-1},\M^{-1,0}]_{q^{-1}}=0\cr [\M^{-2,0},\M^{-2,1}]_q=0\cr}
{}~~~~~~~~~~~~~~~~~~~~~~if~~N=2n+1.
\ee
\end{itemize}
In the case $N=3$, the relations among the generators $\M^{01},\M^{-10},\k^1$
are simply
\be
\cases{[\k^1,L^{01}]_{q^2}=0\cr
       [\k^1,L^{-10}]_{q^{-2}}=0\cr
       [\M^{01},\M^{-10}]_q=q^{-\f 12}\f{1-(\k^1)^{-1}}{q-q^{-1}}.\cr}
\ee

\subsection{Extending $U_q(so(N))$ to $U_q(e^N)$}

{}~~~The `` infinitesimal '' generators $p^i$ of q-translations and the
generator $\La$ of dilatations satisfy the commutation relations
 \cite{fioeu} reported below.
\be
[\k^h,p^i]_{a_{h,i}}=0,~~~~~~~h=1,2,...,n;
\ee
\be
[\M^{01},p^0]=-q^{-1}p^1~~~~~~~~~~~~~[\M^{-1,0}p^0]=p^{-1}~~~~~~~~~~~~~~~~~~
{}~~~~~~~if~~N=2n+1,
\ee
and in all the remaining cases
\be
[\M^{1-m,m},p^i]_{b_{m,i}}=q^{\rho_m}(\delta^i_{-m}-\delta^i_{m-1})p^{i+1}
{}~~~~~~~~~~~[\M^{-m,m-1},p^i]_{b_{m,i}}=q^{\rho_m}(\delta^i_{1-m}-\delta^i_m)
p^{i-1},
\ee
where
\be
a_{m,i}:=q^{2(\delta^i_m-\delta^i_{-m})},~~~~~~~~~~~b_{m,i}:=
(a_{m-1,i})^{\f 12}(a_{m,i})^{-\f 12}.
\ee

The commutation relations of $p$'s among themselves
are those of a quantum space $\rn_q^N$, ${\cal P}_{A~hk}^{~~ij}p^h p^k =0$,
where ${\cal P}_A$ is the projector appearing with negative eigenvalue in the
projector decomposition of the $\hat R$ matrix of $SO_q(N)$ (the
$q$-antisymmetrizer); they
 amount respectively to
\be
p^ip^l=qp^lp^i,~~~~~~~~-l\neq i<l,
{}~~~~~~~~~~~~~~~~
\ee
and
\be
\sum\limits_{l=-j}^jp^lp_l=(p\cdot p)_j(1+q^{-2\rho_j})
\ee
where
\be
(p\cdot p)_j:=\sum\limits_{l=1}^j p^{-l}p_{-l}+\cases{\f{p^0p_0}{1+q^{-1}}
{}~~~~if~~N=2n+1\cr 0~~~~if~~N=2n, \cr}
{}~~~~~~~j=1,...,n;
\ee
consequently
\be
[(p\cdot p)_j,p^l]=0~~~~~~~~~~~~~|l|\le j.
\ee
We see that the algebra $\hat U_q(e^N)$ generated by $\M,\k,p$ is closed.

Finally
\be
[\La,\vec{p}]_{q^{-1}}=0~~~~~~~~~~~~[\La, \k]=0~~~~~~[\La,\M]=0.
\ee
Note that all the commutation relations are $homogeneous$ in $p$.

{\bf Remark} Note that there exists a natural embedding
$\hat U_q(e^N)\hookrightarrow\hat U_q(e^{N+2})$ obtained
by setting equal to zero all the generators of $p^i,\M^{ij},\k^i$ of
$\hat U_q(e^{N+2})$ where either $i$ or $j$ takes the values $\pm (n+1)$.

The q-deformed analogue of the complex conjugation of the algebra of real
translations and
rotations of the real Euclidean space $\rn^N$ can be introduced
whenever $q\in \rn^+$: for such values of $q$ there exists a complex
conjugation
$*$ which
is consistent
with the algebra relations of $U_q(e^N)$, in other words $U_q(e^N)$ equipped
with
$*$ is a $*$-algebra.

The complex conjugation $*$ acts on the Chevalley generators of
$U_q(e^N)$ in the
following way:
\be
(\k^i)^*=\k^i,~~~~~~~~~~~(\M^{1-k,k})^*=q^{-2}\M^{-k,k-1}~~~~~~k\ge 2,
{}~~~~~~~~~~~~~~~~~(\M^{01})^*=q^{-\f 32}\M^{-10}~~~~~if~~N=2n+1,
\ee
\be
(p^i)^*=p^jC_{ji},~~~~~~~~~~~~~~~~~\La^*=\La^{-1};
\ee
$*$ is extended as an algebra antihomomorphism to all of $U_q(e^N)$, i.e.
$(AB)^*=B^*A^*$.
\footnote{Any definition $\La^*=\alpha\La^{-1}$, $\alpha\in\cn$, is
compatible with the algebra relations (1.21); here we will take $\alpha=1$}

\subsection{New $L$ generators of the Euclidean algebra $U_q(e^N)$}

{}~~~The generators $\M$'s presented in the previous subsection do not commute
with $(p\cdot p)_i$.
For representation-theoretical purposes it is convenient to introduce new
generators $L$'s
instead of the $\M $'s by shifting the latter by some functions of the $p$'s,
in
such a way that
$[L,(p\cdot p)_i]=0$. The $L$'s have no classical analogue.

In section 2 we construct the fundamental Hilbert space
representations of $U_q(e^N)$. One can show (Proposition 2)
that for such representations either $(p\cdot p)_i\equiv 0$
identically, $\forall i\ge h$, or all $(p\cdot p)_i$ are strictly
positive definite.
In the former case the algebra reduces to the homogeneous one $U_q(so(N))$,
in the latter case, which we here consider, it follows that we can define
the inverse of $(p\cdot p)_i$.

We define
\be
\cases{
L^{-m,m+1}:=\M^{-m,m+1}+\f{q^{2\rho_{m+1}+2}}{(1-q^2)(p\cdot p)_m}p^{-m}
p^{m+1}\cr
L^{-m-1,m}:=\M^{-m-1,m}+\f{q^{2\rho_{m+1}+1}}
{(1-q^2)(p\cdot p)_m}p^{-m-1}p^m. \cr}
\ee
(similarly for $L^{12}$). Note that this redefinition is possible only
when $q\neq 1$. The basic property of the new generators is the fact that
(compare with relations (1.14),(1.15))
\be
[L^{-m,m+1},p^i]_{b_{i,m}}=0~~~~~~~~~~~~[L^{-m-1,m},p^i]_{b_{i,m}}=0,
\ee
implying
\be
[L^{-m,m+1},(p\cdot p)_i]=0=[L^{-m-1,m},(p\cdot p)_i]~~~~~~~~~~\forall i,m;
\ee
moreover, it is easy to see
that the $L$'s satisfy the same *-conjugation relations as the $\M$'s.

Let us list now the commutation relations satisfied by the $L$'s.
We can define other roots $L$ starting from
simple ones, just in the same way as we did with the $\M$'s,
using relations (1.2)-(1.4) (with the replacement $\M\rightarrow L$).
Simple roots $L$ can be classified into positive and negative ones according to
the
same convention used for the $\M$'s.

\begin{itemize}

\item Let $k\ge h+1$. The commutation relations between positive and
negative simple roots are
\be
[L^{1-m,m},L^{-k,k-1}]_a=0~~~~~~~~~~~~~~~a=\cases{q^{-1}~~~~m\pm 1=k \cr
1~~~~if~~k\neq m, m\pm 1, \cr}
\ee
\be
[L^{1,2},L^{-2,1}]=\f{(p\cdot p)_2p^1p^1}{(1-q^2)[(p\cdot p)_1]^2},
{}~~~~~~~~
[L^{-1,2},L^{-2,-1}]=\f{(p\cdot p)_2p^{-1}p^{-1}}{(1-q^2)[(p\cdot p)_1]^2},~~
{}~~~~~~~~~~~~~if~~~N=2n,
\ee
\be
\cases{[L^{1-m,m},L^{-m,m-1}]_{q^2}=q^{1+2\rho_m}\f{1-\k^{m-1}(\k^m)^{-1}}
{q-q^{-1}}+C_m~~~~~~~~2\le m\le n\cr
[L^{01},L^{-1,0}]_q=q^{-\f 12}\f{1-(\k^1)^{-1}}{q-q^{-1}}+C_1
{}~~~~~~~if~~N=2n+1\cr}
\ee
where
\be
C_1:=\f{q^{\f 12}}{1-q^2}\left[1+q\f{(p\cdot p)_1}{(p\cdot
p)_0}\right]~~~~~~~~~
if~~~N=2n+1
\ee
\be
C_{m+1:}=\f{q^{2\rho_m}}{1-q^2}\left[1-\f{(p\cdot p)_{m-1}(p\cdot p)_{m+1}}
{[(p\cdot p)_m]^2}\right],~~~~~~~~~~~m\ge 1.
\ee

\item The $[\k,L]$ relations, the Serre relations and the $*$-relations
for the $L$ generators are
the same as those of the $\M$ generators.
\end{itemize}

Summing up, the commutations relations among the $L$'s are the same as those
among the $\M$'s, if we add some `` central charges '' ($C_m$).
Let us compare the two sets of generators of $U_q(e^N)$ $\{p,\M,\k\}$ and
$\{p,L,\k\}$.
The $p$'s close a subalgebra, the $\M,\k$'s too; the $L,\k,$'s alone do not,
but
they close
a subalgebra together with the elements $(p\cdot p)_j$ (and $p^{\pm 1}$, in the
case $N=2n$).
However, the action
of the $L,\k$'s on the $p$'s is essentially trivial (they q-commute with the
$p$'s and commute
with the $(p\cdot p)_j$'s), whereas the $M$'s act non-trivially on the $p$'s.

We collect below the whole set of algebra relations characterizing $U_q(e^3)$:
\be
p^{-1}p^0-qp^0p^{-1}=0~~~~~~~~~p^0p^1-qp^1p^0=0~~~~~~~~~~~p^{-1}p^1-p^1p^{-1}-
(q^{\f 12}-q^{-\f 12})p^0p^0=0.
\ee
\be
[k^1,p^i]_a=0,~~~~~~[L^{01},p^i]_{a^{-\f 12}}=0,~~~~~~~[L^{01},p^i]_{a^{-\f
12}}=0~~~~~~~~~
a=\cases{2~~~~if~~i=1 \cr 0~~~~if~~i=0 \cr -2~~~~if~~i=-1}
\ee
\be
[\k^1,L^{01}]_{q^2}=0=[\k^1,L^{-10}]_{q^{-2}}~~~~~~~~~~~~~
[L^{01},L^{-1,0}]_q=q^{-\f 12}\f{(\k^1)^{-1}+\f{(p\cdot p)_1}{(p\cdot
p)_0}}{q^{-1}-q}
\ee
\be
[\La,\M^{ij}]=0=[\La,\k^i]~~~~~~~~~~~~~~~~[p^i,\La]_q=0.
\ee
\subsection{Casimirs of $\hat U_q(e^N)$}

As in the classical case, $\hat U_q(e^N)$ (the subalgebra generated by
$\M,\k,p$
only)
has $n+1-h$ casimirs; their general form mimics the classical one when given in
terms
of the $q$-epsilon tensor and the covariant generators of $U_q(so(N))$
\cite{fioeu}.
The irreps of $\hat U_q(e^N)$  are characterized by the values of the casimirs.
The simplest casimir is the square momentum casimir
\be
\Omega^0 \equiv (p\cdot p)_n.
\ee
Apart from this we write here explicitly only the remaining
`` Pauli-Lubanski '' casimir $\Omega^1$ for
when $N=3,4$. In the limit $q=1$ it is  given respectively by
\be
\Omega^1=\cases{\varepsilon_{ijk}l^{ij}p^k~~~~~~~~if~~N=3 \cr
                w^hw_h,~~~~w_h:=\varepsilon_{hijkl}l^{ij}p^k~~~~~~~~if~~N=4}
\ee
where we have denoted here by $l^{ij}$ the $so(N))$ generator of rotations in
the plane $ij$.
As in the classical case, $\Omega^1$ will vanish on the singlet representation.

\begin{prop} \cite{fioeu}
When $N=3,4$, the Casimirs $\Omega_1$
in terms of $p,L,\k$ generators take respectively the form
\be
\Omega_1=p^0(\k^1)^{-\f 12}-q(q+1)\f{(p\cdot p)_1}{p^0}(\k^1)^{\f 12}
+q^{\f 12}(1-q)(1-q^2)L^{-1,0}L^{0,1}(\k^1)^{\f 12}p^0
\ee
and
$$
\Omega_1=(L^{-2,1}L^{-1,2})(L^{-2,-1}L^{1,2})\k^2(p\cdot p)_1
+\f{q^{-2}}{(q^2-1)^2}
(p\cdot p)_1\{\k^1(L^{-2-1}L^{12})+(\k^1)^{-1}(L^{-21}L^{-12})\}
$$
\be
+\f{q^{-4}(p\cdot p)_1(\k^2)^{-1}}{(q^2-1)^4}\left[1-q^2\k^2
\f{(p\cdot p)_2}{(p\cdot p)_1}\right]^2
-\f{q^{-2}(p\cdot p)_2}{(1-q^2)^2(p\cdot p)_1}[p^{-1}p^{-1}L^{-21}L^{12}
+p^1p^1L^{-2-1}L^{-12}]\k^2.
\ee
\end{prop}

\sect{The fundamental Hilbert space representations of $U_q(e^N)$}

A $*$-representation $\Gamma$ \cite{schmu}
of a $*$-algebra $A$ on a Hilbert ${\cal H}$ space is essentially a
representation of $A$ such that $\Gamma(a^*)=\Gamma(a)^{\dagger}$
($T^{\dagger}$ is the adjoint of $T$) at least on a dense subset of the
Hilbert space ${\cal H}$. In this section we describe the main features of the
fundamental Hilbert space $*$-representations of $U_q(e^N)$ (denoted by
`` irreps '' in the sequel). In particular we focus on the singlet
representation,
which is the one describing a free zero-spin boson on ${\bf R}^N$ in the limit
$q=1$.
For further details and proofs see Ref.
\cite{fioeu}. We find a basis of ${\cal H}$
and show how the generators of $U_q(e^N)$ are to be represented as operators
on the elements of the basis. We don't deal with premature questions
regarding domains of definition of the operators.
The positivity of the scalar product
\be
\cases{~~~~<u|u>\ge 0,~~~~~, \cr <u|u>=0~~~\Leftrightarrow |u>=0, \cr}
{}~~~~~~~~~~~~~~~~~~~~\forall |u>\in{\cal H}
\ee
will be imposed $apriori$ at each step of our construction, and of course will
be essential in determining the structure of the representations.

\subsection{Choice of the observables}

{}~~~Contrary to the classical case, the momenta $p^i$ don't commute with
each-other, therefore cannot be all chosen
as elements a set of commuting observables in order to study Hilbert spaces of
the irreps of $U_q(e^N)$. On the contrary,
among the commuting observables of a complete set characterizing an
irrep we can always take
\be
p_0,(p\cdot p)_1,...,(p\cdot p)_{n-1},(p\cdot p)_n;
\k^1,...,\k^n~~~~~~~~~~~~~~~~~~(p_0\equiv 0~~~if~~N=2n)
\ee
(in fact we can check that they actually make up a complete set for
the `` singlet '' irrep).
It is easy to realize from the commutation relations of $U_q(e^N)$ in the case
$N=2n+1$ that the sign of the eigenvalues of $p_0$ will be the same within each
irrep.
This will mean that to obtain the q-deformed analog of a classical irrep we
have
to sum
two irreps of $U_q(e^N)$ differing only by the sign of $p_0$.

\begin{prop}
\cite{fioeu} There are only the following two alternatives in ${\cal H}$:
\be\cases{
1)~~(p\cdot p)_i\equiv 0~~~~~~~identically~~~~\forall i=h,h+1,...,n; \cr
2)~~(p\cdot p)_i>0~~~~~~~strictly~~~~~\forall i=h,h+1,...,n . \cr}
\ee
\end{prop}
The case 1) corresponds to a `` trivial '' irrep of $U_q(e^N)$,
i.e. to setting $p_i\equiv 0$; the irrep reduces to an
irrep of $U_q(so(N))$, and therefore
won't be considered in the sequel. In the case 2),
as a consequence of the proposition, $[(p\cdot p)_i]^{-1}$ and the operators
$L$ of section 2.3 will be well-defined.

As an introduction to the results of representation theory \cite{fioeu} for
general $N$, we derive
them in the case $N=3$.

\subsection{The case N=3}

{}~~~~The three observables $(p\cdot p)_1,p_0,log_q(\k^1)$ are respectively
the q-deformed analogues of 1) the
square angular momentum; 2) the momentum component along the $x^0$ direction;
3)
the
total angular momentum component along the same direction; of a one-particle
system
in ${\bf R}^3$. In the case q=1 this is a convenient set of observables  for
instance if the
particle is free or subject
to no other force than the one coming from a magnetic field in the $x^0$
direction.

For our derivation we use the algebra relations in formulae (1.32)-(1.35).
Let ${\cal H}$ denote the Hilbert space of an irrep of $U_q(e^3)$.

{}~~

As a first step, we study the representation of the $p$-subalgebra.
We make an ansatz, assuming existence of an eigenspace $\hat{\cal H}\subset
{\cal H}$
of $p_0,(p\cdot p)_1$
\be
(p\cdot p)_1\hat{\cal H}=M^2q^2\hat{\cal H}~~~~~~~~~~~~p_0\hat{\cal
H}=q^2m\hat{\cal H}
\ee
consisting only of normalizable eigenvectors;
$M^2$ is a nonnegative constant with dimensions of a squared mass which we
assume to be positive
(see proposition 2). Then we find that ${\cal H}$ entirely consists of
eigenspaces of
normalizable eigenvectors, too. The remaining generators $L^{01},L^{-10},\k^1$
of
$U_q(e^3)$ will map each of these eigenspace into itself, since they commute
with
$p_0,(p\cdot p)_1$.

Given $any$ vector $|\psi>\in \hat{\cal H}$, according to eq. (1.20),(1.32)
$|\psi_{\pm r}>:=(p^{\pm 1})^r|\psi>$
($r\in\nn,~~l\le n$) will also be an eigenvector of $p_0,(p\cdot p)_1$.
The eigenvalues of $|\psi>,|\psi_{\pm r}>$ will differ by an
integer power of $q$;
the norm of $|\psi_{-r-1}>$ will be given by
\be
<\psi_{-r-1}|\psi_{-r-1}>=<\psi_{-r}|\psi_{-r}>(M^2-q^{-2r-2}m^2)q^2
\ee
If $0<q<1$, there must exist a $r$ such that, $\forall |\psi>\in\hat{\cal H}$
$(p^{-l})^{r+1}|\psi>=0$, otherwise the above
norm would get negative for large $r$. In other words there must
exist an eigenspace of $p_0,(p\cdot p)_1$,
which is annihilated by $p^{-1}$, we call it ${\cal H}_{\vec{0}}$; this also
fixes the
eigenvalue of $(p_0)^2$ up to a sign. Consequently
\be
(p\cdot p)_1{\cal H}_{\vec{0}}=M^2q^2{\cal H}_{\vec{0}}~~~~~~~~~~~~p_0{\cal
H}_{\vec{0}}=
\pm M[1+q^{-1}]^{\f 12}q^2{\cal H}_{\vec{0}}
\ee
If $q>1$ one would find
similarly that there must exist an eigenspace which is annihilated by $p^1$.
Let $\Pg:=(\pi_0,\pi_1)\in \nn\times\zn$;
the subspaces ${\cal H}_{\Pg}:=\La^{\pi_1}(p^1)^{\pi_0}{\cal H}_{\vec{0}}$ are
eigenspaces
of $p_0,(p\cdot p)_1$. Clearly the maps
$p^{\pm 1}:{\cal H}_{\Pg}\rightarrow {\cal H}_{\Pg\pm e_{0}}$ ($e_0\equiv
(1,0)$)
are invertible ($p^{-1}$ is invertible only in the orthogonal complement of its
kernel ${\cal H}_{\vec{0}}$),
the inverse being
$[p^1]^{-1}=\f{q^{-\f 12}}{(p\cdot p)_1-(p\cdot p)_0}p^{-1}$
and $[p^{-1}]^{-1}=\f{q^{-\f 12}}{(p\cdot p)_1-q^{-2}(p\cdot p)_0}p^1$
respectively, as one can easily check using equations (1.32). Using the
irreducibility of ${\cal H}$
and relations (1.32),(2.6) we arrive at the proposition

\begin{prop}
${\cal H}$ can be decomposed into the direct sum
\be
{\cal H}=\bigoplus\limits_{\Pg\in\nn\times\zn}{\cal H}_{\Pg},
{}~~~~~~~~~~~~~~~~
{\cal H}_{\Pg}:=\La^{\pi_1}(p^1)^{\pi_0}{\cal H}_{\vec{0}}
\ee
of orthogonal eigenspaces ${\cal H}_{\Pg}$ of the observables
$p_0,(p\cdot p)_1$,
\be
(p\cdot p)_1{\cal H}_{\Pg}=M^2q^{2(1+\pi_1)}{\cal H}
_{\Pg}~~~~~~~~~~~~p_0{\cal H}_{\Pg}=\pm M[1+q^{-1}]^{\f 12}
q^{\sum\limits_{k=0}^1(1+\pi_k)}{\cal H}_{\Pg}
\ee
\end{prop}
We can visualize the preceding results by thinking of the vectors of ${\cal
H}_{\Pg}$
as functions of $\vec{p}$ with support concentrated on
the circles ${\cal T}^1\subset \rn^3_{\vec{p}}$ drawn in fig. 2. The arrows in
the
figure show the action of the generators $\La,p^{\pm 1}$.

{}~~

As already anticipated, the generators $L^{01},L^{-10},\k^1$ map each ${\cal
H}_{\Pg}$
into itself; as a second step, we study these maps. We can stick to ${\cal
H}_{\vec{0}}$ since,
due to equations (1.33),(2.7), the maps on the other subspaces can be obtained
from these ones
through application to ${\cal H}_{\Pg}$ of powers in the momenta.

The projections of the casimir $\Omega_1$ (1.38)and of relation $(1.34)_3$ onto
${\cal H}_{\vec{0}}$ read
\be
\Omega=L^{-10}L^{01}(\k^1)^{\f 12}+q^{\f 12}\f{(\k^1)^{-\f 12}-(\k)^{\f 12}}
{(q^2-1)(q-1)}.
\ee
\be
[L^{01},L^{-1,0}]_q=q^{\f 12}\f{q^{-1}+(\k^1)^{-1}}{1-q^2}
\ee
Let $\Omega{\cal H}_{\vec{0}}=\omega{\cal H}_{\vec{0}}$,
$\omega\in\rn$, let
$|\psi>\in {\cal H}_{\vec{0}}$ be an eigenvector
of $\k^1$, $\k^1|\psi>=\mu^2|\psi>$, and define
$|\psi_{\pm m}>:=(L^{0\pm 1})^m|\psi>$; $|\psi_{\pm m}>:=(L^{0\pm 1})^|\psi>$
are new eigenvectors
with eigenvalues $\mu^2 q^{\pm 2m}$. The norm of $|\psi_{m+1}>$ reads
\be
<\psi_{m+1}|\psi_{m+1}>=q^{-\f 32}<\psi_m|L^{-10}L^{01}|\psi_m>
=<\psi_m|\psi_m>\{\omega q^{-m}\mu^{-1}+q^{\f 12}\f{1-\mu^{-2}q^{-2m}}{
(q^2-1)(q-1)}\}
\ee
Since $0<q<1$, there must exist a $m$ such that $|\psi_{m+1}>=0$ otherwise the
above
norm would get negative for large $r$. In other words there must
exist a $w\in \rn$ and a highest weight vector $|\vec{0},0>^w\in{\cal
H}_{\vec{0}}$ such that
\be
L^{01}|\vec{0},0>^w=0,~~~~~~~~~~\k^1|\vec{0},0>^w=q^w|\vec{0},0>^w,
{}~~~~~~~~~~~\omega=q^{\f 12}\f{q^{-w}-q^w}{(1-q^2)(1-q)}.
\ee
If we repeat the
same argument with $|\psi_{-m}>:=(L^{-10})^m|\psi>$, we see that its norm
keeps positive for large $m$, hence there exists no lowest weight
vector. On the contrary, if it were $q>1$ there would
exist a lowest weight vector and no highest weight one.

Defining normalized vectors $|\vec{0},j>^w:=N_j(L^{01})^j|\vec{0},0>^w$, we
find
from irreducibility
that they form a basis of ${\cal H}_{\vec{0}}$ and that
$\k^1|\vec{0},j >^w=q^{2j+w}|\vec{0},j >^w$,
$L^{0\pm 1}|\vec{0},j>^w\propto|\vec{0},j\pm 1>$. The coefficient of
proportionality
in the latter relation can be found using formula (2.11) after replacing
$|\psi_m>$ by
$|\vec{0},j >^w$. This completes the study of the structure of ${\cal
H}_{\vec{0}}$.

Now we can easily extend this knowledge to the other spaces ${\cal H}_{\Pg}$.
By
defining
$|\Pg,j >^w$ as the normalized vector proportional to
$\La^{\pi_1}(p^1)^{\pi_0}|\vec{0},j-\pi_0 >^w$, one determines an orthonormal
basis
of ${\cal H}$ consisting of eigenvectors of the commuting observables
$p_0,(p\cdot p)_1,\k^1$. We collect the results in the

\begin{theorem}
A basis ${\cal B}_q^{\w}$ of ${\cal H}^{\w}$ in the case $N=3$ is the set
$\{|\Pg;j,>^w\}$ ($\Pg\in\nn\times\zn$,
$j\in{\cal J}(\Pg):=\{j\in \zn~~|~~j\le \pi_0\}$, with the following
properties.
\be
\cases{ p_0|\Pg;j>^w=\pm M[1+q^{-1}]^{\f 12}
q^{\sum\limits_{k=0}^1(1+\pi_k)}|\Pg;j>^w, \cr
(p\cdot p)_1|\Pg;j>^w=M^2 q^{2(1+\pi_1)}|\Pg;j>^w \cr
\k^1|\Pg;j>^w=q^{2j_1+w}|\Pg;j>^w;\cr}
\ee
Moreover
\be
(\La)^{\pm 1}|\Pg;j>^w=|\Pg\pm\e_1;j>^w,
\ee
\be
p^1|\Pg;j>^w=M[1-q^{2(\pi_0+1)}]^{\f 12}q^{1+\pi_1}|\Pg+\e_0;j+1>^w
\ee
\be
p^{-1}|\Pg;j>^w=M[1-q^{2\pi_0}]^{\f 12}q^{-\rho_l+1+\pi_1}|\Pg-\e_0;j-1>^w,
\ee
(we have set all the arbitrary phase factors equal to 1).
Here $e_0\equiv(1,0)$, $e_1\equiv(0,1)$, $M$ is a constant with
dimensions of a mass, defined modulo integer $q$-powers and characterizing
the irrep together with $w$.  Finally
\be
\cases{
L^{0,1}|\Pg;j>^w= q^{-\f
32-\pi_0}\left[\f{(1+q^{-w-j+\pi_0})(1-q^{-j+\pi_0})}{(q^{-1}-1)(1-q)}\right]^{\
f 12}
|\Pg;j+1>^w \cr
L^{-1,0}|\Pg;j>^w=q^{-\pi_0}\left[\f{(1+q^{1-w-j+\pi_0})(1-q^{1-j+\pi_0})}{(q^{-
1}-1)(1-q)}\right]^{\f 12}
|\Pg;j-1>^w. \cr}
\ee
\end{theorem}
We have appended suffix $w$ to specify that the value of the casimir $\Omega_1$
given
by relation (2.12).

\subsection{The general case}

{}~~~~We introduce a sort of Borel decomposition of $U_q(e^N)$. As a
consequence
of its existence,
Hilbert space representations will be of highest weight type.

{\bf Definition}
We denote by $U_q^{+,N}$ the subalgebra of $U_q(e^N)$ generated by
the positive roots $L$'s and by $p^{-l}$, $l>h(N)$; by $U_q^{-,N}$
the subalgebra generated by the negative roots $L$'s, by $\La^{\pm 1}$,
by $p^l$, $l\ge h(N)$ and, in the case
$N=2n$ only, by $p^{- 1}$. Clearly $U_q(e^N)=U_q^{-,N}\ot U_q^{+,N}$.

\begin{theorem}\cite{fioeu}
The subspace ${\cal H}_G$ of `` highest weight vectors '', i.e.
\be
{\cal H}_G:=\{|\phi>\in {\cal H}~~|~~u|\phi>={\bf 0},~~~~~~~~~~~\forall u\in
U_q^{+,N} \}
\ee
is infinite-dimensional.
A basis of ${\cal H}_G$ is provided by the vectors
$\{(\La)^s|\phi>,~~s\in \zn\}$ and
$\{(\La)^s(p^{\pm 1})^r|\phi>,~~~s\in \zn,~~r\in \nn$\} in the cases
$N=2n+1$ and $N=2n$ respectively; $|\phi>$ is any nontrivial
vector of ${\cal H}_G$.
$|\phi>$ is cyclic in ${\cal H}$ w.r.t. the subalgebra $u_q^{-,N}$.
(In the sequel by `` the highest weight vector '' we will mean
a particular one of these vectors).
The eigenvalues $k^i$ of the operators $\k^i$ are of the type
$k^i=q^{2j_i}\lambda_1$, $j_i\in\zn$, and the constant $\lambda_1$,
$1\ge\lambda_1>q^2$, is a function of the casimirs characterizing the irrep.
\end{theorem}

The existence of highest weight vectors follows when $0<q<1$ from the
requirement of nonnegativity of the scalar product and from the Borel
decomposition given
at the beginning of section 2.3.
The theorem is proved considering first the Hilbert space representations of
the
$p$-subalgebra,
then the Hilbert space representations of the subalgebra of the $L,\k$'s
within each eigenspace of the observables $(p\cdot p)_j$'s; this is possible
because
of formula (1.26).
Contrary to the case of representation theory of
$U_q(so(N))$, in each such eigenspace
there is no lowest weight vector in ${\cal H}_{\vec{0}}$, due
to the presence
of non-vanishing $C_m$'s in the commutation relations (1.29); therefore each
such eigenspace
is infinite-dimensional.

In the sequel we stick to irreps
characterized by $\lambda_1=1,q$ (among which we can find those having
classical analogue).
For this class of irreps we can introduce a vector $\w\in \zn^n$ such that
$\k^i|\phi>=q^{w_i}|\phi>$.
The vector $\w$ depends on the casimirs and together with mass-scale $M$
(defined modulo $q^2$)
completely characterizes an irrep.
We will therefore attach it as a superscript to the symbol ${\cal H}$
and write ${\cal H}^{\w}$.
Now we can formulate the main proposition of this section.

\begin{theorem}
A basis ${\cal B}_q^{\w}$ of ${\cal H}^{\w}$ is the set
$\{|\Pg;\J,\alpha>\}$ ($\Pg\in\nn^{n-h}\times\zn$, $\J\in{\cal J}$,
$\alpha\in A$) with the following properties.
\be
\cases{ p_0|\Pg;\J,\alpha>=\pm M[1+q^{-1}]^{\f 12}
q^{\sum\limits_{k=0}^n(1+\pi_k)}|\Pg;\J,\alpha>,~~~~~~~~~~~
if~~~N=2n\!+\!1; \cr
(p\cdot p)_i|\Pg;\J;\alpha>=M^2 q^{\sum\limits_{k=i}^n
2(1+\pi_k)}|\Pg;\J;\alpha>~~~~~~~~~~~i\ge 1\cr
\k^i|\Pg;\J;\alpha>=q^{2j_i+w_i}|\Pg;\J;\alpha>;\cr}
\ee
Moreover
\be
(\La)^{\pm 1}|\Pg;\J,\alpha>=|\Pg\pm\e_n;\J,\alpha>,
\ee
\be
p^l|\Pg;\J,\alpha>=M[1-q^{2(\pi_{l-1}+1)}]^{\f 12}q^{\sum\limits_{k=l}^n
(1+\pi_k)}|\Pg+\e_{l-1};\J+\y_l,\alpha'>
\ee
\be
p^{-l}|\Pg;\J,\alpha>=M[1-q^{2\pi_{l-1}}]^{\f 12}q^{-\rho_l+
\sum\limits_{k=l}^n(1+\pi_k)}|\Pg-\e_{l-1};\J-\y_l,\alpha'>,
\ee
\be
p^{\pm 1}|\Pg;\J,\alpha>=Mq^{\sum\limits_{k=1}^n(1+\pi_k)}|\Pg;\J\pm y_1,
\alpha'>,~~~~~~~~~~~~~~~if~~~N=2n
\ee
(we have set all the arbitrary phase factors equal to 1).
Here $l>h(N)$,$\e_l\in \nn^{n+1-h},\y_i\in \zn^n$
with $(\e_l)^j=\delta_l^j$, $(\y_i)^j=\delta^j_i$, $M$ is a constant with
dimensions of a mass, defined modulo integer $q$-powers and characterizing
the irrep.  Finally
\be
\cases{
L^{1-m,m}|\Pg;\J,\alpha>=D_m(\Pg,\J,\alpha)|\Pg;\J+y_m-y_{m-1},\alpha'> \cr
L^{-m,m-1}|\Pg;\J,\alpha>=D'_m(\Pg,\J,\alpha)|\Pg;\J-y_m+y_{m-1},\alpha'>
\cr}
\ee
The domain ${\cal J}$ of $\J$ is
\be
{\cal J}:=\{\J\in \zn^n~~|~~j_i\le \pi_{i-1},~~~~i=h+1,h+2,...,n;
{}~~~~~~j_1\in\zn~~~~if~~N=2n\}.
\ee
The coefficients $D_m,D'_m$ and the values of $\alpha'$
depend on the particular irrep under consideration.
$D_m=0$ if $j_i=\pi_{i-1}$.
$A=A(\w,\Pg,\J)$ is a finite set and $\alpha\in A$ are additional
labels identifying the eigenvalues of the observables, if any,
which have to be
added to the ones of formula (100) to get a complete set. $A$ is
trivial (i.e. it has only one element, therefore label $\alpha$ can be omitted)
1) for any value of $\Pg,\J$ if $\w=0$; 2) whenever $j_i=\pi_{i-1}$.
\end{theorem}

{\bf Remarks}.
\begin{itemize}
\item A basis of the subspace ${\cal H}_G$ (see theorem 4)
is $\{|\Pg=s\e_n;\J=\vec{0};\alpha_0>|~~s\in\zn\}$ if $N\!=\!2n\!+\!1$,
$\{|\Pg=s\e_n;\J=r\y_1;\alpha_0>|~~s,r\in\zn\}$ if $N=2n$.

\item As expected the spectra of $(p\cdot p)_i$ are discrete;
they are particularly simple, since they consist only of q-powers.
Note that none of them contains the zero eigenvalue (but the latter is an
accumulation point of the spectra); in particular $(p\cdot p)_n>0$ always,
i.e. `` there is no state in which the nonrelativistic quantum particle
is at rest ''.

\item Theorem 5 summarizes the essential features of the promised
`` q-latticization '' in momentum space $\rn^N_{\vec{p}}$. For each vector
$\Pg$, the equations $(p\cdot p)_l=M^2q^{\sum\limits_{k=l}^n2(1+\pi_k)}$,
$l=h, h+1,...,n$ single out a n-Torus  submanifold ${\cal T}^n$
within $\rn^N_{\vec{p}}$, on which the
vectors with fixed $\Pg$ have support; let us call
${\cal H}_{\Pg}\in {\cal H}$
the linear span of such vectors. The additional specification
of a  vector $\J$, however, selects in ${\cal H}_{\Pg}$ vector(s) having
well-defined angular momentum components $\k^i$, but no well-defined
$p$-angles; in other words the support of each state is $not$ concentrated
on a point of ${\cal T}^n\subset\rn^N_{\vec{p}}$.
For no choice of a complete set of commuting
observables the corresponding eigenvectors would have a point-like support
in $\rn^N_{\vec{p}}$, since no such set can include $N$ functions of the
(non-commuting variables) $p^i$'s. The q-lattice $\{(\Pg,\J),\alpha\}$ has
to be understood in a space where $n+1-h$ dimensions
(corresponding to the first $n+1-h$ observables (2.2))
 are of `` momentum '' type, and the remaining are  of
`` angular momentum '' type. The action of the generators $p,L,\La^{\pm 1}$
on a vector $|\Pg;\J,\alpha>$ can be visualized as a mapping of
the point $(\Pg;\J)$ of the q-lattice into one of its nearest
neighbour points.
\end{itemize}

{\bf Definition} We define the
singlet Irrep as the one characterized by the highest weight
$\w=0$.

With straightforward computations one can verify that the in the singlet
representation
coefficients $D_m,D_m'$ ($m\ge 2$) appearing in formula (2.24)  read
\begin{eqnarray}
D_m(\Pg,\J)=q^{-1-\rho_m-\pi_{m-1}+\pi_{m-2}}\left[\f{(1-q^{j_{m-1}
-\pi_{m-2}-j_m+\pi_{m-1}})
(1-q^{j_{m-1}-\pi_{m-2}-j_m+\pi_{m-1}-2})}{(1-q^2)^2}\right]^{\f 12}
\nonumber\\
D_m(\Pg,\J)=q^{1-\rho_m-\pi_{m-1}+\pi_{m-2}}\left[\f{(1-q^{j_{m-1}-
\pi_{m-2}-j_m+\pi_{m-1}+2})
(1-q^{j_{m-1}-\pi_{m-2}-j_m+\pi_{m-1}})}{(1-q^2)^2}\right]^{\f 12}
\end{eqnarray}
(we have set all arbitrary phase factors equal to 1).

It is easy to verify that by making the tensor product
$(\tilde \Gamma^{\vec{u}},\tilde {\cal H}^{\vec{u}})$ of the singlet
Irrep $(\Gamma^{\vec{0}},{\cal H}^{\vec{0}})$ of $\hat U_q(e^N)$
and an Irrep $(\Gamma^{\vec{u}}_{hom},{\cal H}^{\vec u}_{hom})$
of $U_q^N\equiv U_q(so(N))$ with highest weight $\vec{u}$  we find
a reducible Hilbert space  representation of $\hat U_q(e^N)$ characterized
by the same momentum scale $M$ , as it occurs in the classical case.
To find the Irreps contained in $\tilde \Gamma^{\vec{u}}$ one proceeds
as in the classical Lie algebra representations; namely,
using orthogonality, one determines all the
highest weight vectors contained in ${\cal H}^{\vec{u}}$.

\begin{prop}
\cite{fioeu} Possible highest weights are
of the form $\w\equiv\vec{u}-l\y_1$, $0\le l\le 2u_1$,
if $N=2n+1$,
$\w\equiv\w(l,l'):=\vec{u}-l\cdot sign(u_2-u_1)(\y_2-\y_1)-l'(\y_2+\y_1)$,
$0\le l\le |u_2-u_1|$,  $0\le l'\le u_1+u_2$, if $N=2n$; $\vec{u}$
denote weights of $U_q^N$ . In particular, when $N=3,4$
the sets $\{\w\}$ of weight satisfy the relations $\{w_1\}= \zn$,
$\{\w\}\subset \zn\ot \zn$ respectively.
We have the following tensor product decomposition
\be
\tilde {\Gamma}^{\vec{u}}=\cases{\bigoplus\limits_{l=0}^{2u_1}
\Gamma^{\vec{u}-l\y_1}~~~~~~~~if~~~N=2n+1 \cr
\bigoplus\limits_{0\le l\le |u_2-u_1|; \atop 0\le l'\le u_1+u_2}
\Gamma^{\w(l,l')}~~~~~~~~if~~~N=2n+1 \cr}
\ee
Highest weight vectors can be easily determined from the above described
tensor product construction procedure.
\end{prop}
Note that only the irreps with $w=0,1,2,...$ have  classical analogue.

We give an intuitive picture of the physical content of the spectra of the
observables (2.2) in the singlet representation. The subspace
${\cal H}_i^{\vec{0}}:=\bigoplus\limits_{\{\Pg,~|~\pi_{i-1}=0\}}
{\cal H}_{\Pg}^{\vec{0}}$ is the eigenspace of the observable
$p^{-i} p_{-i}=(p \cdot p)_i-(p \cdot p)_{i-1}$ with the minimum
eigenvalue compatible with a given eigenvalue of $(p \cdot p)_i$,
namely $p^{-i} p_{-i}=M^2q^{\sum\limits_{k=i}^n2(1+\pi_k)}(q^2-1)$;
the latter quantity never vanishes when $q\neq 1$. This means that the there is
always a
`` point zero '' momentum component available in the plane of the
coordinates $i,-i$.
Now let us ask in which `` directions '' of this plane this point zero
momentum component can be pointed.

For the above vhoice of the momenta, the admitted eigenvalues of
$ln_{q^2}(\k^i)$,
i.e. of the angular momentum component in the plane, are $j_i\le 0$
(see formula (2.25)) and show that (except when $N=2n$, $i=1$)
only a `` clockwise '' or `` radial '' orientation are possible.
The anticlockwise is excluded!
If $N=3$, for instance, minimum $p^1p_1$ means that the momentum is
`` almost pointed '' in the $p^0$ direction; $j_1$ represents the
$p^0$-direction component of the (orbital) angular momentum and cannot take
positive values.
This amounts to
sort of a purely `` kinematical '' PT (parity+time-inversion) asymmetry of the
allowed momentum
space (under a PT transformation $\vec{p}$ would remain unchanged, whereas the
angular momentum
component $h$ would change sign). This
is a surprising feature for a lattice theory; in fact, at least
usual equispatiated lattice theories, which are commonly used
nowadays for regularization purposes, cannot have a parity
asimmetry by a well-known no-go-theorem \cite{nie}.
In next section we will see in which sense
in the classical limit $q\rightarrow 1$, however, parity
symmetry is recovered.

\subsection{Configuration space realization}

One can show \cite{fioeu} that the singlet irrep can be realized in
`` $\FR$-configuration space ''. By this we mean that the vectors of
${\cal H}^{\vec{0}}$ can be realized as  elements (`` wave-functions'')
of $\FR$ and the elements of
$U_q(e^N)$ as q-differential operators acting on them. Actually
in Ref. \cite{fioeu}
we give two equivalent `` $\FR$-configuration space '' realizations,
which we call the unbarred and the barred. The scalar product of two vectors
of  ${\cal H}^{\vec{0}}$ is realized as a $\FR$-integral
involving both the barred and unbarred corresponding wavefunctions. These
`` $\FR$-configuration space '' realizations should  be useful for
future functional analysis studies, such as an intrinsic
characterization of ${\cal H}^{\vec{0}}$, questions regarding in
concrete cases the domains of definition of operators representing
some elements of $U_q(e^N)$, etc.
For further deatils we refer the reader to Ref. \cite{fioeu}.

\sect{Classical limit of the singlet irrep}

{}~~~~In order that the representations of $U_q(e^N)$
presented in this work can be considered as
physically realistic, they should describe a system of one free
particle on $\rn^N$ and $U(so(N))$-spin $\w$ in the limit
(understood in some reasonable sense) $q\rightarrow 1$. For simplicity, let
us stick to the case of the singlet irrep $\w=0$, drop the superscript
${\vec{0}}$
and introduce a subscript $q$ on the ket symbols:  we will write $|ket>_q$
instead of $|ket>^{\vec{0}}$.

The commuting observables
\be
p_0,(p\cdot p)_1,...,(p\cdot p)_{n-1},(p\cdot p)_n;
h_1,...,h_n~~~~~~~~~~~~~~~~~~(p_0\equiv 0~~~if~~N=2n+1),
\ee
 (where $h_i:= log_{q^2}(\k^i)$)
make up a complete set both when $q\neq 1$ and $q=1$. We have chosen
$h_i$ instead of $\k^i$ because it is the set of generators
$\{\M^{i,j},h_i,p^i\}$ which has
classical commutation relations in the limit $q\rightarrow 1$.
The eigenvectors $|\Pg;\J>_q$ of the observables (2.2) form the orthonormal
basis
${\cal B}_q$ of theorem 3  for
all $q\in \rn^+-\{1\}$; when $q=1$ the vectors of  ${\cal B}_{q=1}$ are
(orthogonal) distributions,
i.e. elements of the space of functionals
on some space of smooth functions on $\rn^n$, e.g. ${\cal S}(\rn^N)$. We ask
whether
they can be obtained by some sort of limiting procedure when $q\rightarrow 1$.

For each eigenvector $|\Pg,\J>_q$ the
eigenvalues $j_i$ of $h_i:=log_{q^2}(\k^i)$ don't depend on $q$ and
are integers; whereas the eigenvalues
of $(p\cdot p)_i$ (non-uniformly) `` collapse '' to $M^2$ (see formula
$(2.19)_3$):
$$
\lim_{q\rightarrow 1}c_i(q,\Pg)=M^2
{}~~~~~~~~~~where~~~~~(p\cdot
p)_i|\Pg,\J>=:c_i(q,\Pg)|\Pg,\J>~~~~~~i=1,2,...,n,
$$
\be
\lim_{q\rightarrow 1}c_0(q,\Pg)=M
{}~~~~~~~~~~where~~~~~p_0|\Pg,\J>=:c_0(q,\Pg)|\Pg,\J>~~~~~~if~~N=2n+1.
\ee
If we kept $\Pg,\J$ fixed and let $q\rightarrow 1$,  $|\Pg,\J>_q$ would remain
a
normalized
eigenvector and its eigenvalues $c_i$ would go to $M^2$ (independently of
$\Pg$). Consequently,
 the above limit $cannot$ be given a literal sense
${\cal B}_{q=1}=\{\lim_{q\rightarrow 1}|\Pg,\J>_q~,~|\Pg,\J>_q\in {\cal
B}_q\}$.

However, we note that for each fixed $0\le q<1$ and each $\mu_i\in\rn$ there
exist $\Pg$
large enough
such that $c_i(q,\Pg)$ are close to $\mu_i$, and the difference can be made
smaller
and smaller as $q$ approaches 1, because in that limit the point density of the
set $\{q^n\}_{q\in\zn}$ gets
higher and higher around each fixed point on the real
axis. This suggests a more adequate notion of  `` representation limit'', as
given below. First of all,
for each distribution
$|\vec{\mu};\J>_c\in {\cal B}_1$ of the classical representation  defined by
$$
(p\cdot p)_i|\vec{\mu};\J>_c=\mu_i|\vec{\mu};\J>_c,~~~~~~~~~~~~~~
h_i|\vec{\mu};\J>_c=j_i|\vec{\mu};\J>_c~~~~~~~~~~~~~~
<\vec{\mu}',\J'|\vec{\mu},\J>_c=\delta^{\J\J'}\delta(\vec{\mu}-\vec{\mu}')
$$
\be
(and~~~~~~p_0|\vec{\mu};\J>=\mu_0|\vec{\mu};\J>~~~~~~if~~N=2n+1)
\ee
($i=1,...,n$, $\mu_{i+1}\ge \mu_i\ge 0$, $\J\in \zn^n$ and the second $\delta$
is a Dirac's $\delta$),
we can find a vector function
$\tilde{\Pg}(\vec{\mu},q)$
such that
\be
\mu_i=\lim_{q\rightarrow 1}c_i(q,\tilde{\Pg}(\vec{\mu},q)).
\ee
It is easy to see  that this condition is fulfilled e.g. by
\be
\tilde{\pi}_i:=\left[ log_{q^2}(\f{\mu_i}{\mu_{i+1}})\right]~~~~~~~~~~~(and~~~~
\tilde{\pi}_0:=\left[
log_{q^2}(\f{(\mu_0)^2}{\mu_1})\right]~~~~if~~N\!=\!2n\!+\!1)
\ee
($[a]\in\zn$ denotes the integral part of $a\in \rn$, and $\mu_{n+1}\equiv
M^2$).

Therefore we are led to define
\be
\Vert\vec{\mu};\J>_q:=\alpha(q,\vec{\mu})|\tilde{\Pg}(\vec{\mu},q);\J>_q,
\ee
and choose the function $\alpha(q,\vec{\mu})$ in such a way that
\be
\lim\limits_{q\rightarrow
1^-}<\vec{\mu}';\J'\Vert\vec{\mu};\J>_q=\delta(\vec{\mu}'-\vec{\mu})
\delta^{\J\J'};
\ee
the latter limit is in the sense of convergence of
$<\vec{\mu}';\J'\Vert\vec{\mu};\J>_q$ in the space of functionals on smooth
functions of the
two variables $\vec{\mu}',\vec{\mu}$. This is finally the adequate notion of
the
limit we were
looking for. Simbolically
\be
|\vec{\mu};\J>_c=\lim\limits_{q\rightarrow 1^-}\Vert\vec{\mu};\J>_q.
\ee
It is easy to check that a choice of $\alpha$ satisfying relation (3.7) is
\be
\alpha(q,\vec{\mu}):=\prod\limits_{i=1}^n [{\mu_i}(q^{-2}-1)]^{-\f
12}\cdot\cases{1~~~~if~~N=2n \cr
[{\mu_0}\f{(q^{-1}-1)}2]^{-\f 12}~~~~if~~N=2n+1.\cr}
\ee
Let us verify that
$D(\vec{\mu}',\vec{\mu}):=<\vec{\mu}';\J\Vert\vec{\mu};\J>_q$
is
really a $\delta$-convergent functional.

Suppose first that $N=2n$. We consider a smooth function $f(\vec{\mu}')$
and the integral $\int d\vec{\mu}'~D\cdot f(\vec{\mu}')$; we want to show that
its limit is
$f(\vec{\mu})$. We note that
\be
<\tilde{\Pg}(\vec{\mu}',q);\J|\tilde{\Pg}(\vec{\mu},q);\J>_q=\prod_i\chi_{
[q^{2\tilde{\pi}_i(\vec{\mu},q)},q^{2\tilde{\pi}_i(\vec{\mu},q)-2})}(\f{\mu'_i}{
\mu'_{i+1}}),
{}~~~~~~~~~~~~~~\chi_{[a,b)}(z):=\cases{1 ~~~~if~~z\in [a,b)\cr 0~~~~
otherwise.\cr}
\ee
Therefore, setting $z_i:=\f{\mu_i'}{\mu_{i+1}'}$, $i=1,...,n$,
 we find
\be
\prod\limits_{i=1}^n~d\mu_i=\prod\limits_{i=1}^n~dz_i J(\vec{z})
{}~~~~~~~~~J(\vec{z}):=M^N\cdot z_n^{n-1}z_{n-1}^{n-2}...z_2
\ee
and
$$
\int d\vec{\mu}'~D(\vec{\mu}',\vec{\mu})\cdot f(\vec{\mu}')
=\alpha(q,\vec{\mu})M^N
\int d\vec{z}J(\vec{z})\alpha(q,\vec{\mu}'(z))f(\vec{\mu}'(\vec{z}))
\prod\limits_{i=1}^n\chi_{[q^{2\tilde{\pi}_i(\vec{\mu},q)},q^{2\tilde{\pi}_i(\ve
c{\mu},q)-2})}(z_i)
(z_0)
$$
$$
\stackrel{q\rightarrow 1^-}{\ap}|\alpha(q,\vec{\mu})|^2M^Nf(\vec{\mu})
\int d\vec{z}J(\vec{z})
\prod\limits_{i=1}^n\chi_{[q^{2\tilde{\pi}_i(\vec{\mu},q)},q^{2\tilde{\pi}_i(\ve
c{\mu},q)-2})}(z_i)
(z_0)
$$
\be
=\f{(n)_{q^{-2}}!}{n!}
(q^{-2}-1)^nM^{2n}q^{2\sum\limits_{l=1}^nl\tilde\pi_l(\vec{\mu},q)}
|\alpha(q,\vec{\mu})|^2f(\vec{\mu})\stackrel{q\rightarrow 1^-}{\ap}
(q^{-2}-1)^n(\prod\limits_{i=1}^n\mu_i)~|\alpha(q,\vec{\mu})|^2~f(\vec{\mu})
\ee
and the last expression goes to $f(\vec{\mu})$ in the limit $q\rightarrow 1^-$,
due to
relations (3.9). Similarly one proves the result in the case $N=2n+1$.

In the limit $q\rightarrow 1$ the `` PT asymmetry '' in the spectrum
of the observables $\k^i$ noticed at the end of section 2.3 disappears,
`` almost everywhere '' in momentum space.
(Actually, the only points where this does not occur are charaterized by
the condition $(p\cdot p)_{i-1}=(p\cdot p)_i$, namely
they belong to a cylinder in the classical momentum space $\rn^N_{\vec{p}}$;
the latter is a subset of $\rn^N_{\vec{p}}$ of zero measure). In fact,
whenever $\mu_{i-1}<\mu_i$, i.e. $(p\cdot p)_{i-1}<(p\cdot p)_i$,
the range of each $j_i$ (as a function of the
square momenta) becomes the whole set $\zn$ in the limit  $q\rightarrow 1^-$,
since then
$lim_{q\rightarrow 1^-}\tilde\pi_{i-1}(\vec{\mu},q)=\infty$.
The same is true also in the irreps with highest weight $\neq 0$.

\section*{Acknowledgments}
I thank S. Sachse for useful discussions. I am grateful to J.
Wess for his kind hospitality at his Institute and to the Alexander Von
Humboldt-Stiftung for
financial support.

\end{document}